\documentclass{article}

\usepackage{amsmath,amsthm,amsfonts,graphicx,amssymb}
\usepackage{float,xspace,color,array}
\usepackage{pstricks,pst-node,pst-tree}
\newcommand{\perm}{one-stack sortable permutation\xspace}
\newcommand{\perms}{one-stack sortable permutations\xspace}

\newcommand{\Perms}{One-stack sortable permutations\xspace}
\newcommand{\dperm}{decorated one-stack sortable permutation\xspace}
\newcommand{\dperms}{decorated one-stack sortable permutations\xspace}
\newcommand{\N}{\ensuremath{\mathbb{N}}\xspace}
\newcommand{\dis}{\ensuremath{d(\sigma_1,\sigma_2)}}

\newcommand{\T}{{\mathcal T}}

\newcommand{\ens}[1]{\ensuremath{\{#1\}}}
\newcommand{\barre}[1]{\ensuremath{[#1]_a}}

\newtheorem{defn}{Definition}
\newtheorem{thm}{Theorem}
\newtheorem{lem}{Lemma}
\newtheorem{corollary}{Corollary}
\newtheorem{rk}{Remark}
\newtheorem{prop}{Proposition}

\newcommand{\ie}{i.e.}

\begin{document}
%%%%%%%%%%%%%%% TITLE %%%%%%%%%%%%%
\title{Edit distance between unlabeled ordered trees}

%%%%%%%%%%%%%%% NAMES OF THE AUTHORS (WITHOUT ADRESSES)
\author{Anne Micheli \and Dominique Rossin}
%%\email{\{amicheli,rossin\}@liafa.jussieu.fr}
%\address{CNRS, LIAFA, Universit\'e Paris 7, 2 Place
%Jussieu, 75251 PARIS Cedex 05, FRANCE\\Email: \{amicheli,rossin\}@liafa.jussieu.fr}
%\subjclass{05C12,05C05,05A05,05A15}
%\keywords{Edit distance, trees}
%%%%%%%%%% CREATION OF THE TITLE  %%%%%%%%%%%%%

\maketitle

\begin{abstract}
There exists a bijection between one stack sortable permutations --permuta\-tions which avoid the pattern $231$-- and planar trees. We define an edit distance between  permutations which is coherent with the standard edit distance between trees. This one-to-one correspondence yields a polynomial algorithm for the subpermutation problem for $(231)$ avoiding permutations.

Moreover, we obtain the generating function of the edit distance between ordered trees and some special ones.
For the general case we show that the mean edit distance between a planar tree and all other planar trees is at least $n/ln(n)$.

%% Then, we deduce the average distance between those trees both analytically and bijectively. 

Some results can be extended to labeled trees considering colored
Dyck paths or equivalently colored one stack sortable permutations.
\end{abstract}

\section{Introduction}
The edit distance between two trees is the minimal number of edit operations to transform one tree into the other. The edit operations are deletion (edge contraction), insertion of an edge and relabeling of a vertex.

The main problem is to find efficient algorithms to compute this distance between ordered labeled trees. Many algorithms have been proposed \cite{ZS89,klein98}. The basic idea of all these dynamic algorithms arises from the paper of Zhang and Shasha \cite{ZS89}. Further improvements have been made \cite{klein98}. 

Comparing the structure of molecules and finding the preserved ones during a genetic mutation can be seen as an edit distance problem. The application field of this problem is not restricted to biology: in computer vision, objects are represented by their skeletons -which are trees-, and in computer science, edit distance is used to compare structural similarities between XML documents \cite{XML}.

But no combinatorial interpretation has been made of the edit distance between trees. In this article, we introduce \perms \cite{BM00,West90}. These \perms are $(231)$  pattern-avoiding permutations and we show that they are in one-to-one correspondence with ordered trees. 

 Moreover the edit operations can be easily described in terms of \perms. This leads to a purely combinatorial explanation of the edit distance. 

  Some polynomial algorithms are known to compute the edit distance between trees \cite{ZS89}. By our correspondence, we show that computing the greatest common pattern between two $(231)$-avoiding permutations is also polynomial whereas it is NP-complete for general permutations \cite{BBL98}.

\section{Definitions}
\subsection{\Perms}

We describe in this section an encoding for planar trees. We number the edges of the tree by a postfix  traversal and then read the permutation by a prefix traversal. The obtained permutations are called one stack sortable permutations \cite{BM00,West90}. An alternate definition is the following:

\begin{defn}
Let $n \in \N$, a \perm on $\{1\ldots n\}$ is a permutation $\sigma$ such that $\sigma = InJ$ where $I$ and $J$ are \perms on $\{1\ldots p\}$ and $\{p+1\ldots n-1\}$ respectively. Notice that $I$ or $J$ could be empty. 
\end{defn}

Note that in the sequel, permutations are seen as words.

\begin{thm}
\Perms are in one-to-one correspondence with rooted ordered trees. 
%% Notice that the numbering of the edges corresponds to a postfix Depth First Search Traversal (DFS). Then the \perm is obtained by a prefix Depth First Traversal of the labeled tree. (see Figure \ref{fig:T1524376})

\end{thm}

\begin{proof}
Given a tree $T$ with $n$ edges, number the edges by a postfix Depth First Search Traversal (DFS). Read it again by a prefix DFS. It is clear that the obtained permutation is of the form $InJ$. Moreover $I$ corresponds to the encoding by a postfix DFS of the left subtree as shown in Figure \ref{fig:Tperm}. The same goes for $J$ but its numbers are  shifted.

Conversely, take a \perm $\sigma = InJ$. 
\begin{itemize}
\item If $\sigma = k$ then the corresponding tree is a single edge.
\item If $\sigma = InJ$ then the corresponding tree $T_{\sigma}$ is the tree obtained by taking an edge $e = (xy)$ (corresponding to $n$) where $x$ is the root of $T_{\sigma}$. Since $I$ and $J$ are also \perms, we can recursively build the corresponding trees $T_I$ and $T_J$. Put them at each end of the edge $e$, ie $T_I$ is hanging on $x$ such $e$ is the rightmost edge of $x$, and $T_J$ on $y$. 
\end{itemize}
This construction is unique.

\begin{figure}[H]
\begin{center}
\input{Tperm.pstex_t}
\caption{Coding a tree with a \perm.}
\label{fig:Tperm}
\end{center}
\end{figure}
\end{proof}

If $\sigma$ is a \perm, let $\T(\sigma)$ denote the tree associated to $\sigma$. Conversely, if $T$ is a tree, its associated \perm is denoted by $\Theta(T)$.
Moreover, in the sequel, $\sigma_k$ will either denote the $k$-th letter of the word $\sigma$ or the corresponding edge in $\T(\sigma)$.

\begin{defn}
A \emph{subsequence} of a permutation $\sigma = \sigma_1\ldots \sigma_n$ is a word $\sigma'=\sigma_{i_1}\ldots\sigma_{i_k}$ where $i_1,\ldots,i_k$ is an increasing sequence of elements of $\{1,\ldots,n\}$. 

Let $\Phi$ be the bijective mapping of $\{\sigma_{i_1},\sigma_{i_2},\ldots,\sigma_{i_k}\}$ on $\{1,\ldots,k\}$ preserving the order on $\sigma_{i_l}$.

The \emph{normalized} subsequence (pattern) $\hat{\sigma'}$  is equal to $\Phi(\sigma')$.
\end{defn}

\begin{rk}
The \perms are the permutations avoiding the normalized subsequence (pattern) $231$ \cite{Knuth}.
\end{rk}

\subsection{Edit distance}

We briefly recall the definition of the edit distance between trees. Given two trees, the edit distance is the minimal number of operations necessary to transform one into the other. The operations are:
\begin{itemize}
\item {\bf Deletion} : This is the contraction of an edge; two vertices are merged. Only one label is kept.
\item {\bf Insertion} : This is the converse operation of deletion. 
%\item {\bf Relabeling} : This is the change of a label on a vertex.
\end{itemize}

\begin{figure}[H]
\begin{center}
\input{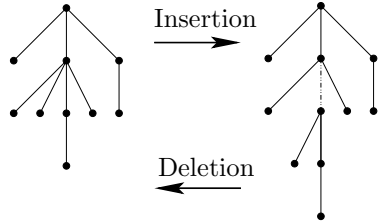}
\caption{Insertion and Deletion operations on a tree.}
\label{fig:InsertDeleteRelabel}
\end{center}
\end{figure}

A cost can be given to each operation. In this article we take $1$ for every cost. 
%For unlabeled trees, the operations considered are {\it deletion} and {\it insertion}.

\section{Distance on \perms}

Since \perms are in one-to-one correspondence with planar trees, we define similar edit operations between \perms and show that these definitions match with edit distance between trees. Moreover, we give a combinatorial interpretation of the distance.

A {\it factor} of a permutation $\sigma=\sigma_1\sigma_2\ldots\sigma_n$ is a {\it factor} of the word $\sigma_1\sigma_2\ldots\sigma_n$ \ie a word  of the form $\sigma_k\sigma_{k+1}\ldots\sigma_{k+l}$.

A factor $f$ is \emph{compact} if it is a permutation of an interval of \N.

A factor $f$ of $\sigma$ is {\it complete} if no non-empty factor $g$ of $\sigma$ verifies both:
\begin{enumerate}
\item  $fg$ is compact where $fg$ is the concatenation of the words $f$ and $g$;
\item  the greatest element of $fg$ is equal to the greatest element of $f$.
\end{enumerate}

Take for example the \perm $\sigma=(1524376)$. The complete factors of $\sigma$ are \ens{1},\ens{15243},\ens{1524376},\ens{5243},\ens{524376},\ens{2},\ens{243},\ens{43},\ens{3},\ens{76},\ens{6}.
\begin{figure}[H]
\begin{center}
$\overset{\mbox{{\tiny\input{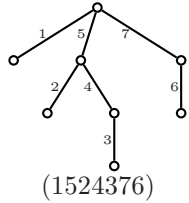}}}}{(1524376)}$
\caption{Tree associated to $\sigma=(1524376)$.}
\label{fig:T1524376}
\end{center}
\end{figure}

A subtree $T'$ of $T$ is a tree such that $T\setminus T'$ is connected.

\begin{lem}
Each compact factor of $\sigma$ are in one-to-one correspondance with:
\begin{itemize}
\item to a subtree
\item to a internal path $P$ in $T = {\mathcal T}(\sigma)$ where each internal vertex of $P$ is of degre $2$ in $T$ and $P$ does not end at a leaf ($P$ can be an internal edge).
\end{itemize}
\label{lem:compact}
\end{lem}

\begin{proof}
First let prove that the subset of edges correpsonding to a compact factor is connected.

Let $\sigma'$ be a compact factor of $\sigma=\Theta(T)$. Let $E_{\sigma'}$ be the set of edges corresponding to $\sigma'$ in $T$. Suppose that $E_{\sigma'}$ is not connected. Let $E_1$ and $E_2$ be two connected components. Let $v$ be the first common ancestor of $E_1$ and $E_2$. Let $P_1$ (resp. $P_2$) be the path starting from $v$ and ending at the first vertex of $E_1$ (resp. $E_2$). Note that we can choose $E_1$ and $E_2$ such that edges of $P_1$ and $P_2$ are not in $E_{\sigma'}$. Suppose that $P_1$ is at the left of $P_2$ (See Figure \ref{fig:compact}).
\begin{figure}[ht]
\begin{center}
\input{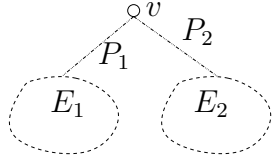}
\label{fig:compact}
\end{center}
\caption{Compact factors are connected components.}
\end{figure}
In the prefix DFS of $T$, edges of $P_2$ are visited between those of $E_1$ and $E_2$. Thus they should appear in $\sigma'$, hence $P_2 = \varnothing$. Thus $v \in E_2$ so that $P_1$ links $E_2$ and $E_1$. In the postfix DFS, the edges of $P_1$ have labels greater than those of $E_1$ and less than $E_2$. If $P_2 \neq \varnothing$, it implies that $\sigma'$ is not compact. Thus $E_{\sigma'}$ is connected.

\begin{figure}[H]
\begin{center}
\input{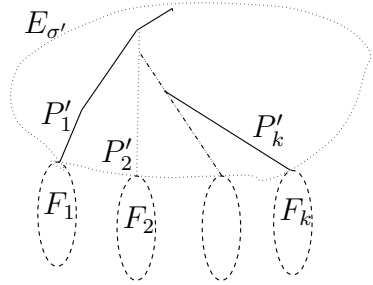}
\caption{Subtree of $T$ induced by $E_{\sigma'}$.}
\label{fig:TTPrime}
\end{center}
\end{figure}

Consider the subtree $T'$ of $T$ induced by $E_{\sigma'}$. It consists of $E_{\sigma'}$ plus all vertices of $T$ that have an ancestor in $E_{\sigma'}$ as shown in
Figure \ref{fig:TTPrime}. 

$E_{\sigma'}$ can be decomposed into edge-disjoint paths $P_i$ thanks to the prefix DFS (See Figure \ref{fig:TTPrime}). $F_i$ is the subtree pending on $P_i$ which can be empty.

The prefix DFS of $T'$ (which is a factor of $\sigma$) gives the associated permutation  $\Theta(P'_1)\Theta(F_1)\Theta(P'_2)$ $\Theta(F_2)\ldots
\Theta(P'_k)\Theta(F_k)$. So $\sigma' = \Theta(P'_1)\Theta(F_1)\Theta(P'_2)\Theta(F_2)\ldots\Theta(P'_k)$, hence $ F_i = \varnothing,\forall i < k$.
\begin{itemize}
\item Suppose $F_k \neq \varnothing$. If $k > 1$, then the edges of $F_k$ are visited after at least one edge of $P'_1$, and before the edges of $P'_k$ in the postfix DFS. Since $\sigma'$ is compact, it implies $k=1$.
\item If $F_k = \varnothing$, $E_{\sigma'}$ is a subtree.
\end{itemize}

The converse is straightforward.

\end{proof}

\begin{prop}
\label{prop:completeTree}
The set of complete factors of $\sigma$ corresponds to the set of subtrees of the associated tree. 
\end{prop}
\begin{proof}
Let $T'$ be a subtree of $T$ and $\sigma = \Theta(T)$. The edges of $T'$ are visited consecutively by the postfix (resp. prefix) DFS of $T$. Thus the sequence of edges of $T'$ is a compact factor $\sigma_k\sigma_{k+1}\ldots\sigma_{k+l}$ of $\sigma$. $\sigma_{k+l+1}$ is an edge which is visited after all edges of $T'$ by the prefix DFS. 
 Thus it is the first time this edge  is visited by the traversal. Hence, its label is greater than those of $T'$. Thus  $\sigma_k\sigma_{k+1}\ldots\sigma_{k+l}$ is complete.
%%%%%%%%%%%%%%%%%%%%%%%%%%%%%%%%%%%%%%%%
%%%%%%%PREUVE EN REEECRITURE %%%%%%%%%%
%%%%%%%%%%%%%%%%%%%%%%%%%%%%%%%%%%%%%%%%

Conversely, let $\sigma'$ be a complete factor. As $\sigma'$ is compact, by Lemma \ref{lem:compact}, it corresponds either to a subtree or to an internal path $P$ with a subtree $F$ hanging on $P$. $\Theta(P)\Theta(F)=\sigma'\Theta(F)$ is also a compact factor of $\sigma$ and it has the same maximum as $\sigma'$ which contradicts the completeness of $\sigma'$.

\end{proof}
\begin{rk}
Let $\sigma$ be a \perm and $\sigma_k = (p(v_k)v_k)$ an edge where $p(v_k)$ denote the parent of $v_k$.
 Let $\sigma'$ be the shortest complete factor of $\sigma$ such that $\sigma' = \sigma_k \sigma_{k+1} \ldots \sigma_{k+l}$ where $\sigma_i = (p(v_i)v_i)$. By previous proposition $\T(\sigma')$ is a subtree of $\T(\sigma)$. The children of $v_k$  are the vertices $v_{k+i}$ such that $i \leq l$ and $\sigma_k > \sigma_{k+i} > \sigma_{k+1},\sigma_{k+2},\ldots,\sigma_{k+i-1}$.
\end{rk}

Let $\sigma = \sigma_1 \ldots \sigma_k$ be a word of $\{1\ldots n\}$ and $a$ be a letter of $\{1\ldots n\}$. We denote by $\barre{\sigma}$ the word $\sigma'_1 \ldots \sigma'_k$ where 
$$\sigma'_i = \begin{cases}
 \sigma_i \text{ if } \sigma_i < a\\
 \sigma_i+1 \text{ otherwise }
\end{cases}
$$

\begin{figure}[t]
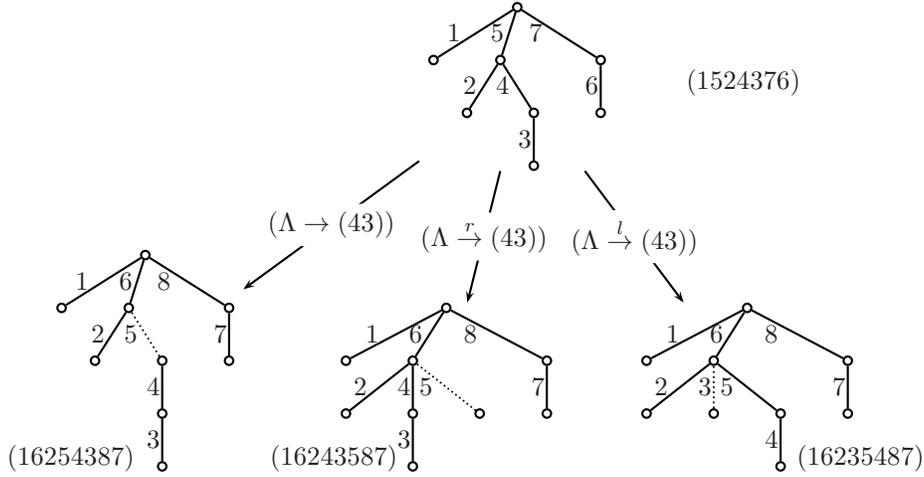

\begin{center}
\begin{pspicture}(0,-0.5)(10,6)
%%\input{Insertion2.pstex_t}
%%\psset{levelsep=100pt}
  %\psset{showbbox=false}
\rput[b](5,4){\rnode{A}{\input{A1524376_0} }}
\rput[b](8,5){\rnode{Alabel}{{$(1524376)$} }}
%\Rnode(6,3){Alabel}{\psframebox*{$(1524376)$ }}
%%  \Tr{\input{A16254387_5}}
 %% \Tr{\input{A16254387_5}}
\rput[b](0,0){\rnode{B}{\input{A16254387_5p}}}
\rput[b](-1,0){\rnode{Blabel}{$(16254387)$}}
\ncline{->}{A}{B}
\ncput*{$(\Lambda \rightarrow (43))$}
\rput[b](4,0){\rnode{C}{\input{A16243587_5p}}}
\rput[b](2.5,0){\rnode{Clabel}{$(16243587)$}}
\ncline{->}{A}{C}
\ncput*{$(\Lambda \overset{r}{\rightarrow} (43))$}
\rput[b](8,0){\rnode{D}{\input{A16235487_3p}}}
\rput[b](9.5,0){\rnode{Dlabel}{$(16235487)$}}
\ncline{->}{A}{D}
\ncput*{$(\Lambda \overset{l}{\rightarrow} (43))$}
\end{pspicture}

\caption{Insertion operations for $f=(43)$.}
\label{fig:insertion}
\end{center}
\end{figure}

 \begin{defn}
%% Let $\sigma = \sigma_1 \ldots \sigma_n$ be a \perm of $S_n$. 
%% Let $k \in \{1\ldots n\}$. We say that:
%% $E = \{\sigma_k \ldots \sigma_{k+l}\}$ belongs to $E_{k,m}$ if its elements are a permutation of an interval of \N and $max\{E\} = m$. Let $E_k = \bigcup_{m \in \{0,\ldots,n\}} \lbrace E, E \in E_{k,m} \text{ and } E \text{ maximal for the inclusion order} \rbrace$                                 
%%%%%%%%%%%%%%%%%%%%%%%%
%% \cup \{\varnothing\}$,
%%%%%%%%%%%%%%%%%%%%%%%%% 

We define two operations on  permutations which map the standard definition on trees (\cite{ZS89}):
\begin{enumerate}
\item{{\it Deletion} : Let $1 \leq k \leq n$. The deletion $(\sigma_k \rightarrow \Lambda)$ is the removal of $\sigma_k$ in a permutation $\sigma$ and the renormalization on $S_{n-1}$ of the result.} We will either talk about the deletion of the edge $\sigma_k$ or the deletion of the vertex $v$ such that $\sigma_k$ is the edge $p(v)v$.
\item{{\it Insertion (see Figure \ref{fig:insertion})} : $(\Lambda \rightarrow \varnothing)$ corresponds to the transformation of the permutation $\sigma=\varnothing$ into $\sigma' = (1)$. If $\sigma \neq \varnothing$,  let $f$ be a complete factor of $\sigma$. Then, $\sigma = ufv$ with $u,v$ factors of $\sigma$.
%%$(\Lambda \rightarrow k)$ with $k \in \{1\ldots n\}$ is one of the following operations:
  \begin{enumerate}
  \item $(\Lambda \rightarrow f)$: The resulting permutation is
  $\sigma' = \barre{u}af\barre{v}$, $a = max\{f\}+1$. This corresponds
  to the insertion of an inner vertex with $\T(f)$ as subtree. 
  \item $(\Lambda \overset{r}{\rightarrow} f)$ : The resulting
  permutation is $\sigma' = \barre{u}fa\barre{v}$, $a =
  max\{f\}+1$. This corresponds to the insertion of a leaf as the
  right sibling of $\T(f)$.
  \item $(\Lambda \overset{l}{\rightarrow} f)$ : The resulting
  permutation is $\sigma' = \barre{u}a\barre{f}\barre{v}$, $a =
  min\{f\}$. This corresponds to the insertion of a leaf as the
  left sibling of $\T(f)$.
  \end{enumerate}
 
}
\end{enumerate}

%% \begin{enumerate}
%% \item{{\it Deletion} : Let $1 \leq k \leq n$. The deletion $(k \rightarrow \Lambda)$ is the removal of $\sigma_k$ in $\sigma$ and the renormalization on $S_{n-1}$ of the result.}
%% \item{{\it Insertion} : $(\Lambda \rightarrow k)$ with $k \in \{1\ldots n+1\}$ defines one of the following operations:
%%   \begin{itemize} 
%%   \item $(\Lambda \rightarrow (k,E))$ with $E \in E_k$
%%     \begin{enumerate}
%%     \item Add $1$ to all elements greater or equal than $i = max\{E\} + 1$. 
%% %%%%%%%%%%%%%%%
%% %%If $E = \varnothing$ then then $i = min\{\sigma_k,\ldots,\sigma_n\}$ : 
%% %%%%%%%%%%%%
%%     \item Insert the element $i $ at the k-th place in $\sigma$.
%%     \end{enumerate}
%%   \item $(\Lambda \rightarrow (k))$ is either one of these two transformations:
%%     \begin{enumerate}
%%     \item    
%%       \begin{enumerate}
	
%%       \item $l = \underset{l' \geq k}{min} \{ l', \sigma_{l'} > \sigma_k\}$. Take $l = n+1$ if $\sigma_{l'} \leq \sigma_k, \forall l' \in \{k,\ldots,n\}$
%%       \item Add $1$ to all elements greater or equal than $i = \sigma_k$.
%%       \item Insert the element $i$ at the l-th place in $\sigma$.
%%       \end{enumerate}
%%     \item 
%%       \begin{enumerate}
%%       \item $i = min\{\sigma_k,\ldots,\sigma_n\}$
%%       \item Add $1$ to all elements greater or equal than $i$. 
%%       \item Insert the element $i $ at the k-th place in $\sigma$.
%%       \end{enumerate}
%%     \end{enumerate}

%%   \end{itemize}
%% }
%% \end{enumerate}

\end{defn}

%% Let $\sigma = 1524376$. Notice that the first choice of values  for $i$ and $l$ in the algorithm - item \ref{noeudInterne} - corresponds to the addition of an edge between two internal nodes. The two other ones - \ref{feuilleG} and \ref{feuilleD} are the addition of an edge between an internal node and a leaf. See figure \ref{fig:insertion} for example.
We study now these operations on the permutation $\sigma=(1524376)$.

\begin{figure}[p]
\includegraphics[width=.9\textwidth]{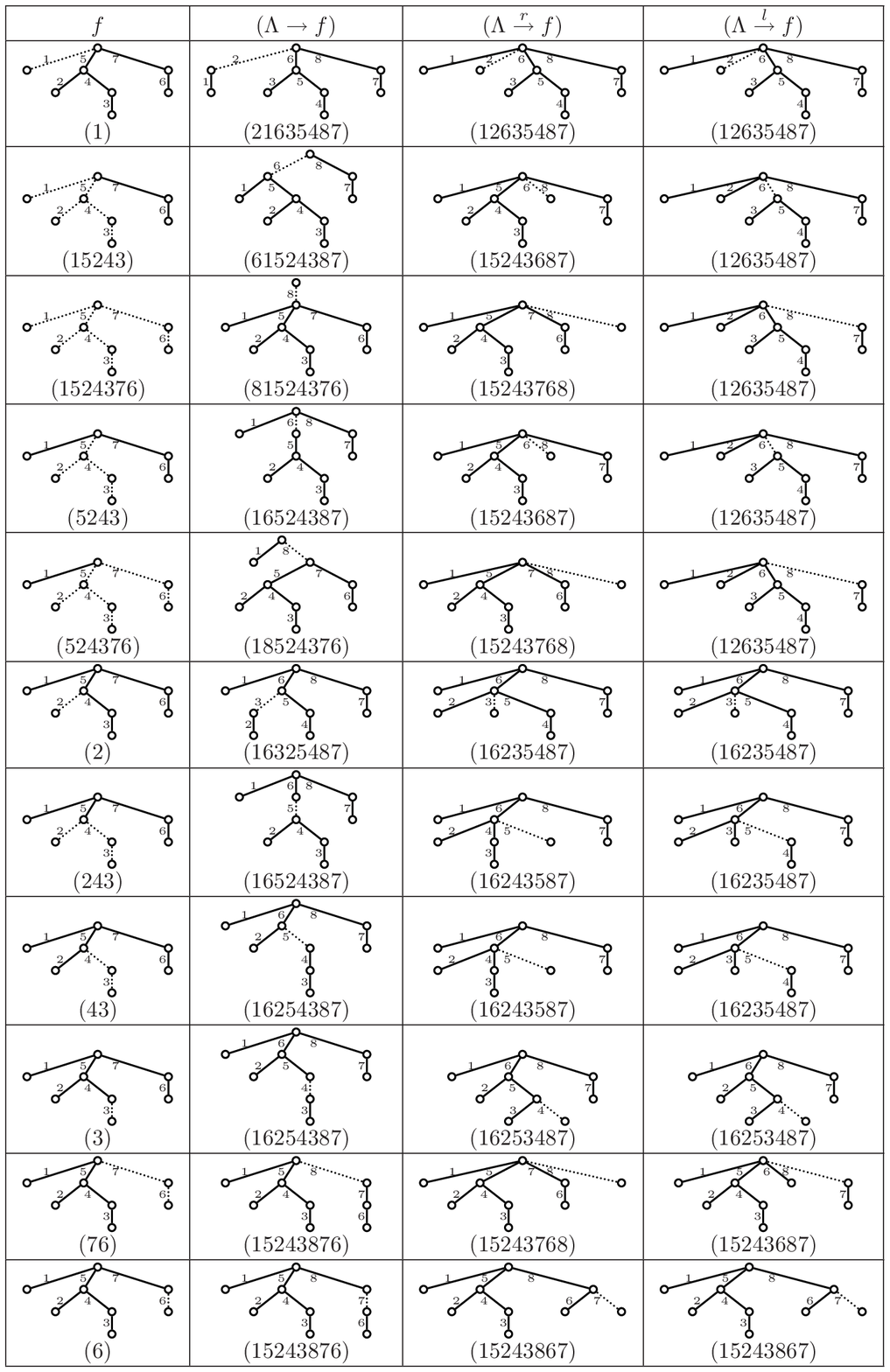}
\caption{Insertion in permutation $\sigma = 1524376$.}
\label{fig:insertions}
\end{figure}

The array of Figure \ref{fig:insertions} gives all the permutations that can be obtained with a single insertion in $\sigma$.

%%\input{A1524376}

%% \begin{figure}[H]
%% $$\begin{array}{|c|c|c|c|}
%% \hline
%% f & (\Lambda \rightarrow f) & (\Lambda \overset{r}{\rightarrow} f) & (\Lambda \overset{l}{\rightarrow} f)\\
%% \hline
%% (1) & 
%% %%\input{A21635487_2}
%% \cellule{21635487}{2} & \cellule{12635487}{2} & \cellule{12635487}{1}\\
%% \hline
%% (15243) & \cellule{61524387}{6} & \cellule{15243687}{6}  & \cellule{12635487}{1}\\
%% \hline
%% (1524376) & \cellule{81524376}{8} & \cellule{15243768}{8} & \cellule{12635487}{1}\\
%% \hline
%% (5243) & \cellule{16524387}{6} & \cellule{15243687}{6} & \cellule{12635487}{2}\\
%% \hline 
%% (524376) & \cellule{18524376}{8} & \cellule{15243768}{8} & \cellule{12635487}{2}\\
%% \hline
%% (2) & \cellule{16325487}{3} & \cellule{16235487}{3} & \cellule{16235487}{2}\\
%% \hline
%% (243) & \cellule{16524387}{5} & \cellule{16243587}{5} & \cellule{16235487}{2} \\
%% \hline
%% (43) & \cellule{16254387}{5} & \cellule{16243587}{5} & \cellule{16235487}{3}\\
%% \hline
%% (3) & \cellule{16254387}{4} & \cellule{16253487}{4} & \cellule{16253487}{3} \\
%% \hline
%% (76) & \cellule{15243876}{8} & \cellule{15243768}{8} & \cellule{15243687}{6}\\
%% \hline
%% (6) & \cellule{15243876}{7} & \cellule{15243867}{7} & \cellule{15243867}{6}\\
%% \hline
%% \end{array}$$
%% \caption{Insertion in permutation $\sigma = 1524376$}
%% \end{figure}

We prove now that the operations (deletion and insertion) defined on \perms are in fact internal operators for \perms. Moreover, these operators define an edit distance between permutations coherent with  the usual edit distance between trees. 

\begin{lem}
The Deletion/Insertion algorithm yields a \perm.
\end{lem}

\begin{proof}

\begin{itemize}
\item Deletion : The proof is straightforward considering the one-to-one correspondence with trees and \perms. Consider a tree labeled by a depth first traversal. Deleting the edge $i$ from this tree changes all labels greater than $i$ by subtracting~$1$. 
\item Insertion : 
Let $\sigma$ be a \perm and $f$ be a complete factor of $\sigma = ufv$. By Proposition \ref{prop:completeTree}, $f$ corresponds to a subtree of $\T(\sigma)$. 
\begin{enumerate}
\item $(\Lambda \rightarrow f)$: Let $T = \T(\sigma)$ and
$(e_1,e_2,\ldots,e_n)$ be the edges of $T$ ordered by a prefix DFS of
the tree. Note that $\sigma = \alpha_T(e_1)\alpha_T(e_2)\ldots \alpha_T(e_n)$ where $\alpha(i)$ is the label of the edge $i$ in $T$.

Let $T'$ be the tree obtained by the insertion of an internal vertex $v$ ($a = (p(v)v)$) at the
root vertex of the subtree ${\mathcal T}(f)$.  Moreover ${\mathcal T}(f)$ is a subtree hanging on
$v$. Let $\sigma'' = \Theta(T')$.  A prefix traversal of $T'$ orders
the edges of $T'$ as follows:
$(e_1,e_2,\ldots,e_l,a,e_{l+1},\ldots,e_n)$.

  Since $\sigma''$ is obtained by a prefix traversal, $\sigma'' =
u'af'v'$. Since the edges of $f$ appear before $a$ in the postfix DFS,
$f'= f$. The edge $a$ in a
postfix DFS appears just after $f$. Thus its label is
$max\{f\}+1$. All the edges visited after $f$ in $T$ (and so after $a$
in $T'$) by the postfix DFS have their labels increased by $1$. Thus
$\sigma'' = \barre{u}af\barre{v} = \sigma'$.

%% Let $T = \T(\sigma)$ and $T'$ be the
%% tree obtained after the insertion of an edge $a$ at the root vertex of
%% the subtree $f$. Let $\sigma'' = \Theta(T')$. Since $\sigma''$ is
%% obtained by a prefix traversal,  $\sigma'' = u'af'v'$. $f'=f$ since
%% the edges of $f$ appear before $a$ by a postfix DFS. We look at the
%% edges (not the labels) visited by the prefix DFS. The edges of $T'$
%% are in the same order than those of $T$.
%%  except the insertion of the edge
%% $a$.
%%  This edge is inserted before the factor $f$. The edge $a$ in a
%% postfix DFS appears just after $f$. Thus its label is
%% $max\{f\}+1$. All the edges visited after $f$ in $T$ (and so after $a$
%% in $T'$) in the postfix DFS have their labels increased by $1$. Thus
%% $\sigma'' = \barre{u}af\barre{v} = \sigma'$.

\item $(\Lambda \overset{l}{\rightarrow}  f),(\Lambda \overset{r}{\rightarrow}  f)$ : The same arguments as for $(\Lambda \rightarrow  f)$ hold.
\end{enumerate} 
\end{itemize}
\end{proof}

\begin{prop}
Insertion and deletion are inverse operations.
  \label{prop:operationInverse}
\end{prop}

\begin{proof}
  There are two different kinds of deletions in a tree $T$.
  \begin{enumerate}
  \item Deletion of an inner vertex $v$. Consider the subtree $T'$ of $T$ hanging on $v$. It corresponds to a complete factor $f$ in $\sigma = \Theta(T)$. This contraction corresponds to the inverse operation of $(\Lambda \rightarrow f)$.
  \item Deletion of  a leaf. There are three different cases:
    \begin{itemize}
    \item Deletion of a vertex with no sibling. This is the same as deleting the parent of this vertex which is an inner vertex except if the tree is reduced to a single edge.
    \item Otherwise, this vertex has either:
      \begin{itemize}
      \item A left sibling $v'$. Consider the subtree hanging at $v'$ (including $p(v')v'$). It corresponds to the factor $f$. The inverse operation is $(\Lambda \overset{r}{\rightarrow} f)$
      \item A right sibling $v'$. Consider the subtree hanging at $v'$ (including $p(v')v'$). It corresponds to the factor $f$. The inverse operation is $(\Lambda \overset{l}{\rightarrow} f)$
      \end{itemize}
    \end{itemize}
  \end{enumerate}
  
\end{proof}

\begin{defn}
The distance between two \perms $\sigma_1$ and $\sigma_2$ is the minimal number of operations -deletion or insertion - to transform $\sigma_1$ into $\sigma_2$.
\end{defn}

For example let $\sigma_1 = 31264587$ and $\sigma_2 = 1524376$. We want to transform $\sigma_1$ into $\sigma_2$.

\begin{center}
\begin{tabular}{c|cc}
\parbox[t]{5.5cm}{
\begin{itemize}
\item $31264587 \xrightarrow{(1 \rightarrow \Lambda)} 2153476$
\item $2153476 \xrightarrow{(1 \rightarrow \Lambda)} 142365$
\item $142365 \xrightarrow{(\Lambda \rightarrow 3) }1524376$
\end{itemize}
} & \parbox[t]{5cm}{
%%  \begin{figure}[H]
    $\overset{\mbox{{\tiny\def\dedge{\ncline[linestyle=dotted,dotsep=1pt]}
 \psset{levelsep=10pt,labelsep=0pt,tnpos=l,radius=2pt}
\pstree{\TC
}{
\pstree{\TC
\tlput{$3$}}{
\TC[edge=\dedge]
\tlput{$1$}
\TC
\tlput{$2$}
}
\pstree{\TC
\tlput{$6$}}{
\TC
\tlput{$4$}
\TC
\tlput{$5$}
}
\pstree{\TC
\tlput{$8$}}{
\TC
\tlput{$7$}
}
}
}}}{(31264587)}$\\
    $\overset{\mbox{{\tiny\def\dedge{\ncline[linestyle=dotted,dotsep=1pt]}
 \psset{levelsep=10pt,labelsep=0pt,tnpos=l,radius=2pt}
\pstree{\TC[edge=\dedge]
}{
\TC
\tlput{$1$}
\pstree{\TC
\tlput{$4$}}{
\TC
\tlput{$2$}
\TC
\tlput{$3$}
}
\pstree{\TC
\tlput{$6$}}{
\TC
\tlput{$5$}
}
}
}}}{(142365)}$
%%\end{figure}
} & 
\parbox[t]{5cm}{
%%  \begin{figure}
    $\overset{\mbox{{\tiny\def\dedge{\ncline[linestyle=dotted,dotsep=1pt]}
 \psset{levelsep=10pt,labelsep=0pt,tnpos=l,radius=2pt}
\pstree{\TC
}{
\pstree{\TC
\tlput{$2$}}{
\TC[edge=\dedge]
\tlput{$1$}
}
\pstree{\TC
\tlput{$5$}}{
\TC
\tlput{$3$}
\TC
\tlput{$4$}
}
\pstree{\TC
\tlput{$7$}}{
\TC
\tlput{$6$}
}
}
}}}{(2153476)}$\\
    $\overset{\mbox{{\tiny\def\dedge{\ncline[linestyle=dotted,dotsep=1pt]}
 \psset{levelsep=10pt,labelsep=0pt,tnpos=l,radius=2pt}
\pstree{\TC
}{
\TC
\tlput{$1$}
\pstree{\TC
\tlput{$5$}}{
\TC
\tlput{$2$}
\pstree{\TC[edge=\dedge]
\tlput{$4$}}{
\TC
\tlput{$3$}
}
}
\pstree{\TC
\tlput{$7$}}{
\TC
\tlput{$6$}
}
}
}}}{(1524376)}$
%%  \end{figure}
}
%%  & 
%% \parbox[t]{3cm}{
%% %%  \begin{figure}
%%     \input{A142365}
%% %%  \end{figure}
%% } & 
%% \parbox[t]{3cm}{
%% %%  \begin{figure}
%%     \input{A1524376_4}
%% %% \end{figure}
%% }
\\
\end{tabular}
\end{center}

\begin{thm}
The edit distance between ordered trees is the distance between the associated \perms. 
\end{thm}
\begin{proof}
This is a consequence of Proposition \ref{prop:operationInverse}.
\end{proof}

\begin{thm}
The edit distance between \perms $\sigma_1$ and $\sigma_2$ is equal to 
$$|\sigma_1| + |\sigma_2| - 2 |u|$$
where $u$ is a largest normalized subsequence (pattern) of $\sigma_1$ and $\sigma_2$.
\end{thm}

\begin{proof}
The edit distance $\dis$ between $\sigma_1$ and $\sigma_2$ is given by the minimal number of insertions and deletions.  If $t_1$ is an insertion and $t_{2}$ is a deletion then there exist a deletion $t'_1$ and an insertion $t'_{2}$ such that $t_1t_{2}(\sigma)=t'_{1}t'_2(\sigma)$. Note that $t'_1$ and $t'_2$ depend on the \perm $\sigma$.

Considering the sequence of edit operations, there exists a sequence made of deletions then insertions that transforms $\sigma_1$ into $\sigma_2$. We denote this sequence by $D_1\ldots D_lO_1\ldots O_k$, $l+k = \dis$.

Consider the \perm $\sigma' = D_1\ldots D_l(\sigma_1)$. Take $u =
\sigma'$. $u$ is a normalized subsequence of $\sigma_1$ because
deleting an edge from a \perm yields a normalized subsequence of the
original \perm.  $u$ is also a normalized subsequence of $\sigma_2$
because inserting an edge in a \perm $s$ yields a \perm $s'$ and $s$
is a normalized subsequence of $s'$.

Conversely, take $u$ as a maximal normalized subsequence of $\sigma_1$
and $\sigma_2$. It is straightforward to find $|\sigma_1| - |u|$
operations of deletions such that those deletions transform
$\sigma_1$ into $u$. The same goes for $\sigma_2$ and $u$.

\end{proof}

\begin{corollary}
Finding the greatest common pattern between two \perms is polynomial.
\end{corollary}

In \cite{BBL98}, they proved that finding the greatest common pattern between two permutations is NP-complete. We prove here that the problem becomes polynomial when restricting to \perms, ie $(132)$ or $(231)$-avoiding permutations. In fact, the algorithm of Zhang and  Shasha \cite{ZS89} on trees solves the problem on \perms because the algorithm outputs not only the distance but also the greatest common subtree.

%%%%%%%%%%%%%%%%%%%%%%%%%%%%%%%%%%%%%%%%%%%%%%%%%%%%%%%%%%%%%%%%%%%%
%%%%%%%%%%%%%%%%%%%%%%%%%%%%%%%%%%%%%%%%%%%%%%%%%%%%%%%%%%%%%%%%%%%%
%% Section distance mmoyenne min
%%%%%%%%%%%%%%%%%%%%%%%%%%%%%%%%%%%%%%%%%%%%%%%%%%%%%%%%%%%%%%%%%%%%
%%%%%%%%%%%%%%%%%%%%%%%%%%%%%%%%%%%%%%%%%%%%%%%%%%%%%%%%%%%%%%%%%%%%
\section{Lower bounds on average edit distance}

In this section we study the average edit distance between a given planar tree $T$ with $n$ vertices and all other
 planar trees with $n$ vertices. We show that this average distance is lower bounded by $\frac{n}{ln(n)}$.

\begin{lem}
Let $T$ be a planar tree with $n$ vertices. There are at most $n-1$ different deletions and $3n^3$  insertions allowed in $T$.
\end{lem}

\begin{proof}
The number of deletions is upper bounded by the number of edges \ie $n-1$.

The number of insertions is bounded by $3$ times the number of subtrees (or complete factor of the corresponding permutation). The number of subtrees of $T$ rooted at vertex $v$ is bounded by $d(v)^2$ where $d(v)$ denotes the degree of vertex $v$. Thus the total number of subtrees is bounded by $\sum_v d(v)^2$. 
\end{proof}

\begin{thm}
Let $T_0$ be a tree with $n$ vertices. The proportion of planar trees with $n$ vertices at distance at most ${\mathcal O}(n/ln(n))$ tends to $0$.

The average distance between $T_0$ and the set of planar trees is lower bounded by $n/ln(n)$.
\end{thm}

\begin{proof}
Let $T_0$ be a planar tree. Let $A_k = \{ T \in {\mathcal T}_n , \text{dist}(T_0,T) \leq 2k \}$. Note that $A_0 = \{ T_0 \}$. A tree $T_k \in A_k$ is obtained from $T_0$ by $l \leq k$ deletions then $l$ insertions. Thus $|A_k| < (n-1)^k (n^3)^k < n^{4k}$. But the number of planar trees $C_n \equiv \frac{4^n}{n \sqrt{\pi n}}$. So that the proportion of planar trees at distance at most ${\mathcal O}(n/ln(n))$ tends to $0$. 

Hence the average distance is lower bounded by $n/ln(n)$.
\end{proof}

\section{Generating functions}

\begin{figure}[H]
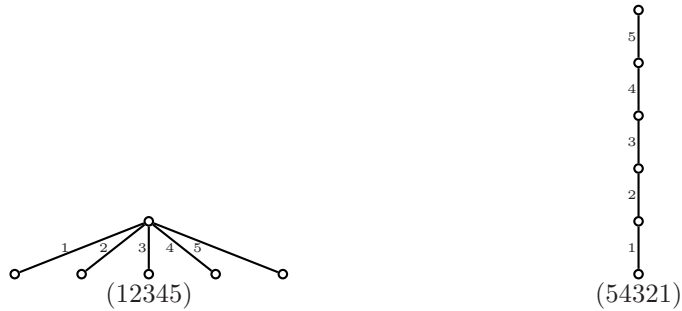

\begin{center}
   \psset{levelsep=20pt,labelsep=0pt,tnpos=l,radius=2pt}

 $$ \overset{\mbox{{\tiny\input{A12345_0}}}}{(12345)}  \hspace*{4cm} \overset{\mbox{{\tiny\input{A54321_0}}}}{(54321)} $$
\caption{Some canonical trees.}
\label{fig:peigne}
\end{center}
\end{figure}

Using the combinatorial interpretation of the distance, we compute the
generating functions of the edit distance between planar trees with
$n$ edges and some special ones as
shown in Figure \ref{fig:peigne}. Moreover, we deduce the average
distances from the generating functions.

\subsection{Generating function of the edit distance between \perms and $Id~=~1~2~\ldots~n$}

We denote by $S_1(t,q)$ the generating function of \perms where $t$
counts the size of the permutation and $q$ the edit distance
between \perms and $Id$.  This is the distance between a tree and the
trivial one which is made of $n$ edges and of height $1$.

%%%%%%%%%%%%%%%%%%%%%%%%%%%%%%%%%%%%%%%%%%%%%%%%%%%%%%%%%%%%%%%%%%%%%%
%%%%%%%%%%%%%%%%%%%%%%%%%%%%%%%%%%%%%%%%%%%%%%%%%%%%%%%%%%%%%%%%%%%%%%
%%%%%%%%%%%%%%%%%%%%%%%%%%%%%%%%%%%%%%%%%%%%%%%%%%%%%%%%%%%%%%%%%%%%%%

\paragraph{{\bf Tree interpretation of the largest increasing subsequence}}

\begin{prop}
The length of a largest increasing subsequence of a \perm is the number of leaves of the associated tree.
\label{prop:croissant}
\end{prop}

\begin{proof}
Let $T$ be a planar rooted tree and $\sigma$ the associated \perm. We call a leaf-edge an edge incident to a leaf. 

\begin{enumerate}
\item The subsequence of $\sigma$ made of the leaf-edges is increasing because the order in which the leaf-edges are visited by a prefix traversal is the same than by a postfix traversal.

\item Suppose that we take an increasing subsequence $\sigma'$ of $\sigma$. This subsequence is in one-to-one correspondence with some edges in the tree. Suppose that there is an internal one $\gamma = (p(\nu)\nu)$. Then, by the postordering of the edges, each edge $(p(v)v)$ such that $\nu = p(v)$ has a smaller label and appears in $\sigma$ after the edge $\gamma$. Thus, none of these edges are in $\sigma'$. Moreover, there is at least one leaf edge belonging to the subtree $T_\gamma$ hanging on $\nu$. Replace edge $\gamma$ by a leaf of $T_\gamma$. The prefix traversal ensures that the obtained subsequence is an increasing one.

\end{enumerate}
\end{proof}

\begin{prop}
The number of rooted planar trees with $n$ edges and $k$ leaves is equal to the number of rooted planar trees with $n$ edges and $n+1-k$ leaves. 
\label{prop:symetrie}
\end{prop}

\begin{proof}
This is a direct consequence of the symmetry of the Narayana numbers $\frac{1}{n}\binom{n}{k}\binom{n}{k-1} $ which count the number of planar trees with $n$ edges and $k$ leaves. 

\end{proof}

%%%%%%%%%%%%%%%%%%%%%%%%%%%%%%%%%%%%%%%%%%%%%%%%%%%%%%%%%%%%%%%%%%%%%%
%%%%%%%%%%%%%%%%%%%%%%%%%%%%%%%%%%%%%%%%%%%%%%%%%%%%%%%%%%%%%%%%%%%%%%
%%%%%%%%%%%%%%%%%%%%%%%%%%%%%%%%%%%%%%%%%%%%%%%%%%%%%%%%%%%%%%%%%%%%%%

\paragraph{{\bf Generating function}}

We now compute the generating function $I(t,p)$ of \perms of size $t$ and largest increasing subsequence of size $p$. 
\begin{itemize}
\item $[I(t,p)]_0 = 1$
\item $[I(t,p)]_1 = p$
\item $[I(t,p)]_2 = (p+p^2)$
\end{itemize}
\begin{eqnarray}
[I(t,p)]_{n} = p [I(t,p)]_{n-1} + \sum_{i=0}^{n-2}[I(t,p)]_{i}[I(t,p)]_{n-1-i}
\label{equ:rec}
\end{eqnarray}
This formula comes from the decomposition of a \perm $\sigma$ into $InJ$ with $n \geq 1$. The largest increasing subsequence of $\sigma$ is the union of the largest one of $I$ and the largest one of $J$ unless $J$ is empty - in this case, the largest subsequence is the largest one for $In$ -.

From this formula we deduce: 
%%the associated equality on generating functions:
\begin{eqnarray}
I(t,p) = 1 + (p-1)tI(t,p) + t I^2(t,p)
\label{equ:genI}
\end{eqnarray}
\begin{itemize}
\item $1$ comes from the case $n=0$ in the equation (\ref{equ:rec}). 
\item $ptI(t,p)$ comes from $p [I(t,p)]_{n-1}$~.
\end{itemize}

It follows from equation (\ref{equ:genI}):
\begin{eqnarray}
I(t,p) = \frac{1+(1-p)t-\sqrt{(p-1)^2t^2-2(p+1)t+1}}{2t}
\end{eqnarray}
Let $\tilde{S_1}(t,q)$ be the generating function of the difference between the lengths of the \perm and the largest increasing subsequence in it. 
\begin{itemize}
\item $[\tilde{S_1}(t,q)]_{0} = -1$
\item $[\tilde{S_1}(t,q)]_{1} = 0$
\item $[\tilde{S_1}(t,q)]_{2} = q$
\end{itemize}

\begin{lem}
\begin{eqnarray}
I(t,p) = 1+p+p \tilde{S_1}(t,p)
\end{eqnarray}
\end{lem}
\begin{proof}
\begin{eqnarray*}
I(t,p) &= &\sum_{\tau \geq 1} \sum_{\alpha=1}^{\tau} [I(t,p)]_{\tau,\alpha} t^{\tau} p^{\alpha} +1\\
& = & \sum_{\tau \geq 1} \sum_{\beta=1}^{\tau} [I(t,p)]_{\tau,\tau+1-\beta} t^{\tau} p^{\tau+1-\beta} +1\\
& = & \sum_{\tau \geq 1} \sum_{\beta=0}^{\tau} [I(t,p)]_{\tau,\tau+1-\beta} t^{\tau} p^{\tau+1-\beta} +1 \\
& = & 1 + p (\tilde{S_1}(t,p)+1)
\end{eqnarray*}
The end of the proof is straightforward using Proposition \ref{prop:symetrie}.
\end{proof}

\begin{thm}
\begin{eqnarray*}
S_1(t,q) &= &\tilde{S_1}(t,q^2)\\
& = & \frac{1+(q^2-1)t-\sqrt{(q^2-1)^2t^2-2(q^2+1)t+1}}{2tq^2}
\end{eqnarray*}
\end{thm}

%% \subsubsection{Average distance}

%% The average distance $\delta$ can be obtained from the generating function $S_1(t,q)$ in the following way:

%% \begin{itemize}
%% \item $\left. F(t) = \frac{\partial S_1(t,q)}{\partial q} \right|_{q=1}$
%% \item $\delta = \frac{[F(t)]_n}{C(n)}$ where $C(n)$ is the $n$-th Catalan number.
%% \end{itemize}

%% This easy computation yields $\delta = n-1$ but a direct combinatorial interpretation proves this result in a more comprehensive way.

%% \paragraph{{\bf Average edit distance}}

%% \begin{thm}
%% The average edit distance between rooted planar trees with $n$ edges and $n,n-1,\ldots,2,1$ is $n-1$.
%% \end{thm}

%% \begin{proof}
%% This is a direct consequence of Propositions \ref{prop:croissant} and \ref{prop:symetrie}. Another proof can be found in \cite{DHW03,Rei02}. In \cite{Rei02} the result is more general. Thus we provide here a simpler proof for this special case.
%% \end{proof}

\subsubsection{Average distance}

\begin{thm}
The average edit distance between rooted planar trees with $n$ edges and $Id$ is $n-1$.
\end{thm}

\begin{proof}
\begin{enumerate}
\item The average distance $\delta$ can be obtained from the generating function $S_1(t,q)$ in the following way:

\begin{itemize}
\item $\left. F(t) = \frac{\partial S_1(t,q)}{\partial q} \right|_{q=1}$
\item $\delta = \frac{[F(t)]_n}{C(n)}$ where $C(n)$ is the $n$-th Catalan number.
\end{itemize}

This easy computation yields $\delta = n-1$ but a direct combinatorial interpretation proves this result in a more comprehensive way.
\item This is a direct consequence of Propositions \ref{prop:croissant} and \ref{prop:symetrie}. Another proof can be found in \cite{DHW03,Rei02}. In \cite{Rei02} the result is more general. Thus we provide here a simpler proof for this special case.
\end{enumerate}
\end{proof}

\subsection{Generating function of the edit distance between \perms and $n (n-1) \ldots 1$}

This is the distance between a tree and the trivial one which is made
of $n$ edges and is of height $n$.  It is equivalent to finding the
largest decreasing subsequence in the \perm.

We compute the generating function $D(x,y,z)$ of trees with respect to the number of edges $x$, the height of the tree $y$ and the number of leaves $z$ at maximal depth.

\begin{prop}
\begin{eqnarray}
D(x,y,z) &=&  yD(x,y,\frac{1}{1-xz}) - yD(x,y,1) + \frac{xyz}{1-xz}
\end{eqnarray}
\label{prop:equD}
\end{prop}

\begin{proof}
\begin{eqnarray*}
[D(x,y,z)]_{i,j,k} &= &\sum_{l = 1}^{i-j+1} \binom{l+k-1}{k} [D(x,y,z)]_{i-k,j-1,l} \text{ if } j > 1\\
%%\end{eqnarray}
%%\begin{eqnarray}
{[}D(x,y,z){]}_{i,1,k} &= &\delta_{i,k}
\end{eqnarray*}
The coefficient $[D(x,y,z)]_{i,j,k}$ is equal to the number of ways to
add $k$ leaves at depth $j$ to any tree with $i-k$ edges, depth $j-1$
and $l$ leaves at depth $j-1$. $\binom{l+k-1}{k}$ is the number of
ways to add $k$ leaves to $l$ leaves at depth $j$.

\begin{eqnarray*}
D(x,y,z) &= &\sum_{i \geq 1} \sum_{j \geq 1} \sum_{k \geq 1} d_{i,j,k} x^i y^j z^k \\
&=& \sum_{i \geq 1} \sum_{j \geq 2} \sum_{k \geq 1} \sum_{l \geq 1} \binom{k+l-1}{k} d_{i-k,j-1,l} x^i y^j z^k  + y \sum_{i \geq 1} (xz)^i\\
&=& \sum_{i \geq 1} \sum_{j \geq 2} \sum_{k \geq 1} \sum_{l \geq 1} (-1)^k \binom{-l}{k} d_{i-k,j-1,l} x^i y^j z^k  + y \sum_{i \geq 1} (xz)^i\\
&=& \sum_{i \geq 1} \sum_{j \geq 2} \sum_{k \geq 1} \sum_{l \geq 1} (-1)^k \binom{-l}{k} d_{i,j-1,l} x^i y^j (xz)^k  + y \sum_{i \geq 1} (xz)^i
\end{eqnarray*}
Using 
$$(x+a)^{-n} = \sum_{k = 0}^{\infty} \binom{-n}{k} x^k a^{-n-k}$$
\begin{eqnarray*}
D(x,y,z) &=& \sum_{i \geq 1} \sum_{j \geq 2} \sum_{l \geq 1}  ((1-zx)^{-l}-1) d_{i,j-1,l} x^i y^j   + y \sum_{i \geq 1} (xz)^i\\
&=& yD(x,y,\frac{1}{1-xz}) - yD(x,y,1) + \frac{xyz}{1-xz}
\end{eqnarray*}
\end{proof}

Let $S_2(x,y)$ be the generating function with respect to the length $n$ of the \perm and the edit distance between this \perm and $n(n-1)(n-2)\ldots 1$.
Then, $S_2(x,y) = D(xy^2,\frac{1}{y^2},1)$.
%%F(x,y,1)$ where $F(x,y,z)$ is defined by the following functionnal equation:

%% \begin{eqnarray*}
%% F(x,y,z) &= & \frac{1}{y^2} F(x,y,\frac{1}{1-xz}) -\frac{1}{y^2} F(x,y,1) +\frac{xz}{1-xy^2z}
%% \end{eqnarray*}

%% This comes from Proposition \ref{prop:equD} and the fact that $F(x,y,z) = D(xy^2,\frac{1}{y^2},z)$.

In \cite{BKR72,RV96}, they give a solution for $D(x,y,1)$ in terms of a continued fraction. 
$$D(x,y,1) = \sum D_k(y) x^k, D_k(y) = \cfrac{1}{k\begin{cases}1-\cfrac{y}{1-\cfrac{y}{1-\ldots}}
\end{cases}}$$

This yields the solution for $S_2$.
$$S_2(x,y) = \sum y^{2k}D_k(\frac{1}{y^2}) x^k$$

The first terms of $S_2$ are given by:
$$S_2(x,y) = x+{x}^{2}{y}^{2}+{x}^{2}+{x}^{3}{y}^{4}+3\,{x}^{3}{y}^{2}+{x}^{3}+{x}^{4}{y}^{6}+7\,{x}^{4}{y}^{4}+5\,{x}^{4}{y}^{2}+{x}^{4}$$
\begin{center}
\begin{tabular}{|m{0.5cm}m{0.5cm}|m{4cm}m{1cm}|m{2cm}m{1cm}|}
\hline
\psset{labelsep=0pt,levelsep=20pt,radius=2pt}
%\pstree[levelsep=20pt]{\TC}{\TC} 
$\overset{\mbox{{\tiny\def\dedge{\ncline[linestyle=dotted,dotsep=1pt]}
 %%\psset{levelsep=10pt,labelsep=0pt,tnpos=l,radius=2pt}
\pstree{\TC[edge=\dedge]
}{
\TC
\tlput{$1$}
}
}}}{(1)}$
 & $x$ &&&&\\ &&&&&\\
\hline
\psset{labelsep=0pt,levelsep=20pt,radius=2pt}
%\pstree[levelsep=20pt]{\TC}{\pstree[levelsep=20pt]{\TC}{\TC}} 
$\overset{\mbox{{\tiny\def\dedge{\ncline[linestyle=dotted,dotsep=1pt]}
 %%\psset{levelsep=10pt,labelsep=0pt,tnpos=l,radius=2pt}
\pstree{\TC[edge=\dedge]
}{
\pstree{\TC
\tlput{$2$}}{
\TC
\tlput{$1$}
}
}
}}}{(21)}$
& $x^2$ & \psset{labelsep=0pt,levelsep=20pt,radius=2pt}
%\pstree[levelsep=20pt]{\TC}{\TC\TC} 
$\overset{\mbox{{\tiny\def\dedge{\ncline[linestyle=dotted,dotsep=1pt]}
 %%\psset{levelsep=10pt,labelsep=0pt,tnpos=l,radius=2pt}
\pstree{\TC[edge=\dedge]
}{
\TC
\tlput{$1$}
\TC
\tlput{$2$}
}
}}}{(12)}$
& $x^2y^2$ &&\\&&&&&\\
\hline
\psset{labelsep=0pt,levelsep=20pt,radius=2pt}
%\pstree[levelsep=20pt]{\TC}{\pstree[levelsep=20pt]{\TC}{\pstree[levelsep=20pt]{\TC}{\TC}}}  
$\overset{\mbox{{\tiny\def\dedge{\ncline[linestyle=dotted,dotsep=1pt]}
 %%\psset{levelsep=10pt,labelsep=0pt,tnpos=l,radius=2pt}
\pstree{\TC[edge=\dedge]
}{
\pstree{\TC
\tlput{$3$}}{
\pstree{\TC
\tlput{$2$}}{
\TC
\tlput{$1$}
}
}
}
}}}{(321)}$
& $x^3$ & \psset{labelsep=0pt,levelsep=20pt,radius=2pt}
%\pstree[levelsep=20pt]{\TC}{\pstree[levelsep=20pt]{\TC}{\TC\TC}}
$\overset{\mbox{{\tiny\def\dedge{\ncline[linestyle=dotted,dotsep=1pt]}
 %%\psset{levelsep=10pt,labelsep=0pt,tnpos=l,radius=2pt}
\pstree{\TC[edge=\dedge]
}{
\pstree{\TC
\tlput{$3$}}{
\TC
\tlput{$1$}
\TC
\tlput{$2$}
}
}
}}}{(312)}$
% \pstree[levelsep=20pt]{\TC}{\TC\pstree[levelsep=20pt]{\TC}{\TC}} 
%\pstree[levelsep=20pt]{\TC}{\pstree[levelsep=20pt]{\TC}{\TC}\TC} 
$\overset{\mbox{{\tiny\def\dedge{\ncline[linestyle=dotted,dotsep=1pt]}
 %%\psset{levelsep=10pt,labelsep=0pt,tnpos=l,radius=2pt}
\pstree{\TC[edge=\dedge]
}{
\TC
\tlput{$1$}
\pstree{\TC
\tlput{$3$}}{
\TC
\tlput{$2$}
}
}
}}}{(132)}$
$\overset{\mbox{{\tiny\def\dedge{\ncline[linestyle=dotted,dotsep=1pt]}
 %%\psset{levelsep=10pt,labelsep=0pt,tnpos=l,radius=2pt}
\pstree{\TC[edge=\dedge]
}{
\pstree{\TC
\tlput{$2$}}{
\TC
\tlput{$1$}
}
\TC
\tlput{$3$}
}
}}}{(213)}$
& $3x^3y^2$ & \psset{labelsep=0pt,levelsep=20pt,radius=2pt}
$\overset{\mbox{{\tiny\def\dedge{\ncline[linestyle=dotted,dotsep=1pt]}
 %%\psset{levelsep=10pt,labelsep=0pt,tnpos=l,radius=2pt}
\pstree{\TC[edge=\dedge]
}{
\TC
\tlput{$1$}
\TC
\tlput{$2$}
\TC
\tlput{$3$}
}
}}}{(123)}$
%\pstree[levelsep=20pt]{\TC}{\TC\TC\TC}
& $x^3y^4$\\&&&&&\\
\hline
\end{tabular}
\end{center}
\medskip
\paragraph{{\bf Average edit distance}}

In \cite{BKR72}, they determine analytically the average height of a planar tree with $n$ edges which is ${\sqrt{\pi n}}-\frac{1}{2}$. Thus, the average edit distance is $2 (n-\sqrt{\pi n}+\frac{1}{2}) \equiv 2n$.

\section{Conclusion}
In section 2.2, we define the edit operations to be insertion and deletion. Indeed we omitted a third one, the relabeling operation. Instead of working with unlabeled trees, we study trees whose vertices are labeled and the relabeling operation consists in changing the label of a vertex.

The general case where the trees are labeled and the different edit operations have different costs can be obtained in a similar way. Define a \dperm as a \perm where each number is indexed by a letter; $1_e5_a2_a4_b3_d7_b6_c$ represents the following tree:

\begin{figure}[H]
\begin{center}
\includegraphics[height=3cm]{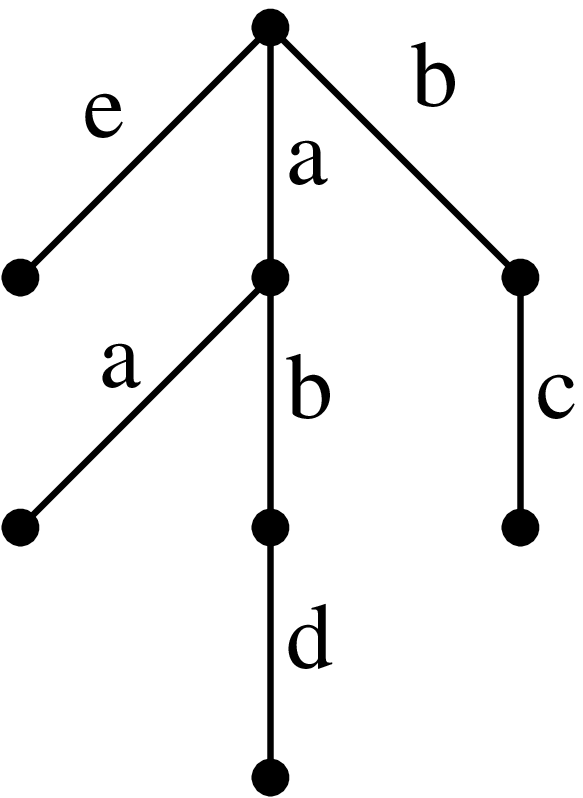}
\end{center}
\end{figure}

The operations on \dperms are almost the same as before and the relabeling operation consists in changing one letter. $c_i,c_d,c_r$ are respectively the insert, delete and relabeling unitary costs. There exists only a difference for the insertion of a new free edge. In the unlabeled case, we did not take into account the insertion of a leaf with no sibling. Thus we define a fourth insertion operation as:
\begin{itemize}
\item $(\Lambda \overset{1}{\rightarrow} i)$ where $i$ is a complete factor of size $1$ of the permutation $\sigma = u i v$. $\sigma' = \barre{u} \barre{i} a \barre{v}$ where $a = i$.
\end{itemize}

Let $\sigma_1$ and $\sigma_2$ be two \dperms with the same underlying permutation. The label distance $d(\sigma_1,\sigma_2)$ is equal to the string distance between both labeled words.

Let $T_1$ and $T_2$ be two \dperms. We denote by a subpermutation $\sigma$ of $T_1$ and $T_2$ a normalized subpermutation without label. $\Sigma_{T_1}$ is the set of all sub-\dperms of $T_1$ which underlying permutation is $\sigma$.

The relabeling distance between $T_1$ and $T_2$ with respect to $\sigma$ is:
$$d_{\sigma}(T_1,T_2) = min\{c_r d(\alpha,\beta), \forall \alpha \in \Sigma_{T_1}, \beta \in \Sigma_{T_2}\}$$

The distance between these two \dperms $T_1$ and $T_2$ is given by $min \{c_{i}(|T_1|-|\sigma|) + c_{d}(|T_2|-|\sigma|) + d_{\sigma}(T_1,T_2), \sigma \text{ normalized subpermutation of } T_1,T_2\}$

\end{document}